\documentclass[12pt]{article}

\usepackage{graphicx}
\usepackage{amsmath}
\usepackage{amsthm}

\newtheorem{theorem}{Theorem}[section]
\newtheorem{definition}{Definition}
\newtheorem{lemma}[theorem]{Lemma}

\begin{document}
	
\textbf{One-sided $S$-systems of operations and character-free proof of Frobenius theorem}

\medskip

\textsc{Kuznetsov Eugene}

\section*{Introduction}

Frobenius groups are one of the classical groups in permutation
group theory. The study of these groups was begun in the Frobenius'
article \cite{Frob02} in the beginning of 20th century and continued
by M.Hall \cite{Hall62}, H.Wielandt \cite{Wiel64} etc. Using the character group theory, Frobenius proved
that \textit{there exists an invariant regular subgroup consisting of all fixed-point-free permutations and the identity permutation in a finite Frobenius group} \cite{Frob02}. This is now known as the \textbf{Frobenius theorem}. A proof of this theorem which does not rely on the character group theory has not been known until now.

\begin{definition}\label{Def1.8}\cite{Hall62,Wiel64} 
The transitive irregular 	permutation group $G$ acting on a set $E$ is called a \textbf{Frobenius group}, if $St_{ab}(G)=<\mathbf{id}>$ for any $a,b\in E,\quad a\neq b$.
\end{definition}

Let $G$ be an arbitrary Frobenius group of permutations on some set
$E$. We use the following notation:

$H_{a}=St_{a}(G)$ is the stabilizer of the element $a\in E$ in the group $G$,

$0$ and $1$ are two distinct elements of the set $E$,
\begin{equation*}
	E^{\ast} = \{0\} \cup \{h(1)\;\vert \;h \in H_{0} = St_{0} (G)\} \subseteq E.
\end{equation*}

In this work, the author continues investigations \cite{Kuzn96} and  gives a character-free proof of the Frobenius theorem. The new proof is based on some notions and results from the theory of ternary operations, the theory
of orthogonal binary operations, the theory of transversals in groups
and the theory of quasigroups and loops. These objects
and theories are described in \cite{Belous65,Belous71,Kuzn941,Kuzn952}.
Proof consists of two sections.

In Section 1, a partial ternary operation $(x,a,y)$ is introduced
and its properties are proved in Lemma \ref{Lemma=00003D00002021}.
In particular, it is proved that a binary operation $A_{a}=(x,a,y)$
is an idempotent quasigroup for any $a\in E^{\ast}\setminus\{0,1\}$.

In the beginning of Section 2 it is proved that the collection of binary
operations $A_{a}(x,y)=(x,a,y)$ is a right $S$-system \cite{Belous65}
of idempotent quasigroups, when $a\neq0,1$ (Lemma \ref{lemma=00002032}).
Then it is proved that for any $a,b\in E^{\ast}$, $a\neq b$, operations
$A_{a}(x,y)=(x,a,y)$ and $A_{b}(x,y)=(x,b,y)$ are orthogonal 
(Lemma \ref{lemma=00002034}). In Lemma \ref{lem34}, two collections
of functions $\varphi_{b,a}$ and $\psi_{b,a,d}$ on the set $E$
are introduced, and it is proved that all of them are permutations.
In Lemma \ref{lem9.8} some properties of these permutations are shown.
Also it is proved that any finite Frobenius group contain a transitive
set of permutations $T$, consisting only fixed-point-free permutations
and identity permutation. Lemma \ref{lemma=00002038} shows that the set $T$ contains all fixed-point-free permutations of the group $G$.

Lemma \ref{lemma=00002039} is the final lemma in this proof of the Frobenius
theorem. It shows that the set $T$ is a left loop transversal
in the group $G$ to its subgroup $H_{0}=St_{0}(G)$. Then the author demonstrated that set $T$ is closed in relation to multiplication and
inverse transformation in the group $G$, i.e. $T$ is a subgroup
in the group $G$. Finally, the author proves that the subgroup $T$ is
an invariant subgroup in the group $G$.

\section{Construction of a ternary operation over a Frobenius group}

\begin{lemma} \label{lemma=00003D00002020}
For any Frobenius group $G$ of degree $n$, the following statements are true: 
\begin{enumerate}
\item There \ exists \ a \ set \ $P=\{\sigma_{x}\}_{x\in E}$ \ of \ $n$ \ permutations \
such \ that $\sigma_{x}(0)=x\quad\forall x\in E$, and $\sigma_{0}=\mathbf{id}$. 
\item For  $\forall a\in E^{\ast}$ and $x\in E$, there exists a unique
permutation $\alpha\in G$ such that 
\begin{equation}
\alpha\in H_{x}=St_{x}(G),\quad\quad h(1)=(\sigma_{x}^{-1}\alpha\sigma_{x})(1)=a.\label{eq1}
\end{equation}
\end{enumerate}
\end{lemma}

\begin{proof}
1. It is obvious because any Frobenius group is a transitive group
on the set $E$.

2. For any $\alpha\in H_{x}$ we have 
\[
h=(\sigma_{x}^{-1}\alpha\sigma_{x})\in(\sigma_{x}^{-1}H_{x}\sigma_{x})=H_{0}
\]
The subgroup $H_{0}$ of the Frobenius group $G$ is sharply transitive on
the set $E^{\ast}\setminus\{0\}$, so $\forall a\in E^{\ast}\setminus\{0\}$
there exists a unique permutation $h_{a}\in H_{0}$ such that $h_{a}(1)=a$.
Then there exists a unique permutation $\alpha_{a}=\sigma_{x}h_{a}\sigma_{x}^{-1}$
satisfying the conditions (\ref{eq1}).
\end{proof}

Let us define an operation $\left(x,a,y\right)$ as follows:
\begin{equation}
\begin{array}{c}
(x,0,y)\overset{def}{=}x,\quad(x,1,y)\overset{def}{=}y,\\
\forall a\in E^{\ast},\;a\neq0,\;1:\quad\quad(x,a,y)\overset{def}{=}\alpha_{a}(y),\\
\text{where \ }\alpha_{a}\in H_{x},\quad h_{a}(1)=(\sigma_{x}^{-1}\alpha_{a}\sigma_{x})(1)=a.
\end{array}\label{eq6.8}
\end{equation}

\begin{lemma}\label{Lemma=00003D00002021} The following statements are true: 
\begin{enumerate}
\item the operation $\left(x,a,y\right)$ is defined correctly; 
\item $(0,a,1)=a,\quad(x,a,x)=x;$ 
\item $\forall a,b\in E^{\ast}:\ (x,a,(x,b,y))=(x,c,y)$ \ for some  $c=c(a,b)\in E^{\ast};$ 
\item if $n$ is finite, then any binary operation $A_{a}=(x,a,y)$ is a quasigroup for any $a\in E^{\ast}\setminus\{0,1\}$. 
\end{enumerate}
\end{lemma}

\begin{proof}
1. See Lemma \ref{lemma=00003D00002020}.

2. We have 
\[
(0,a,1)=u\quad\Leftrightarrow\quad\left\{ \begin{array}{c}
u=\alpha_{a}(1)\\
\alpha_{a}\in H_{0}\\
(\sigma_{0}^{-1}\alpha_{a}\sigma_{0})(1)=a
\end{array}\right.\quad\Rightarrow\quad\left\{ \begin{array}{c}
u=\alpha_{a}(1)\\
\sigma_{0}=\mathbf{id}\\
\alpha_{a}(1)=a
\end{array}\right.\quad\Rightarrow\quad u=\alpha_{a}(1)=a,
\]
i.e. $(0,a,1)=a$. 
\[
(x,a,x)=u\quad\Leftrightarrow\quad\left\{ \begin{array}{c}
u=\alpha_{a}(x)\\
\alpha_{a}\in H_{x}\\
(\sigma_{x}^{-1}\alpha_{a}\sigma_{x})(1)=a
\end{array}\right.\quad\Rightarrow\quad u=\alpha_{a}(x)=x,
\]
i.e. $(x,a,x)=x$.

3. If $a=0$, then we get 
\[
(x,a,(x,b,y))=(x,0,(x,b,y))=x=(x,0,y)\quad\Rightarrow\quad c=c(0,b)=0.
\]

If $a=1$, then we get 
\[
(x,a,(x,b,y))=(x,1,(x,b,y))=(x,b,y)\quad\Rightarrow\quad c=c(1,b)=b.
\]

If $b=0$, then we get 
\[
(x,a,(x,b,y))=(x,a,(x,0,y))=(x,a,x)=x=(x,0,y)\quad\Rightarrow\quad c=c(a,0)=0.
\]

If $b=1$, then we get 
\[
(x,a,(x,b,y))=(x,a,(x,1,y))=(x,a,y)\quad\Rightarrow\quad c=c(a,1)=a.
\]

Let $a,b\neq 0,1$. Then we get 
\[
(x,a,(x,b,y))=u\quad\Leftrightarrow\quad\left\{ \begin{array}{c}
(x,a,v)=u\\
(x,b,y)=v
\end{array}\right.\quad\Leftrightarrow
\]
\[
\Leftrightarrow\quad\left\{ \begin{array}{c}
\alpha_{a}\in H_{x},\quad\alpha_{a}(v)=u\\
(\sigma_{x}^{-1}\alpha_{a}\sigma_{x})(1)=a\\
\alpha_{b}\in H_{x},\quad\alpha_{b}(y)=v\\
(\sigma_{x}^{-1}\alpha_{b}\sigma_{x})(1)=b
\end{array}\right.\quad\Leftrightarrow\quad\left\{ \begin{array}{c}
\alpha_{a},\alpha_{b}\in H_{x},\\
(\sigma_{x}^{-1}\alpha_{a}\sigma_{x})(1)=a\\
\alpha_{a}(\alpha_{b}(y))=v\\
(\sigma_{x}^{-1}\alpha_{b}\sigma_{x})(1)=b
\end{array}\right.
\]
Then for permutation $\alpha=(\alpha_{a}\alpha_{b})$ we get 
\[
\left\{ \begin{array}{c}
u=\alpha(y)=(\alpha_{a}\alpha_{b})(y)\\
\alpha(x)=(\alpha_{a}\alpha_{b})(x)=\alpha_{a}(x)=x\quad\Rightarrow\quad\alpha\in H_{x}\\
\sigma_{x}^{-1}\alpha\sigma_{x}=\sigma_{x}^{-1}(\alpha_{a}\alpha_{b})\sigma_{x}=(\sigma_{x}^{-1}\alpha_{a}\sigma_{x})(\sigma_{x}^{-1}\alpha_{b}\sigma_{x})=h_{a}h_{b}=h_{c}\\
c=h_{c}(1)=h_{a}(h_{b}(1))=h_{a}(b)
\end{array}\right.
\]
where $h_{a}=(\sigma_{x}^{-1}\alpha_{a}\sigma_{x})\in H_{0}$ and $h_{b}=(\sigma_{x}^{-1}\alpha_{b}\sigma_{x})\in H_{0}$. So we obtain 
\begin{equation}
(x,a,(x,b,y))=(\alpha_{a}\alpha_{b})(y)=\alpha(y)=\alpha_{c}(y)=(x,c,y)=(x,h_{a}(b),y).\label{eq=00003D0000203_11}
\end{equation}

4. Let $a\in E^{\ast}\setminus\{0,1\}$. In order to prove that the binary operation $A_{a}=(x,a,y)$ is a quasigroup, it is necessary to show that both equations $(b,a,y)=c$ and $(x,a,b)=c$ have a unique solution in the set $E$ for any fixed parameters $b,c\in E$, $a\in E^{\ast}\setminus\{0,1\}$.

a) Let us have an equation 
\begin{equation}
(b,a,y)=c,\label{eq=00003D0000203}
\end{equation}
where $b,c\in E$, $a\in E^{\ast}\setminus\{0,1\}$. Then 
\begin{equation}
\left\{ \begin{array}{c}
c=\alpha_{a}(y)\\
\alpha_{a}\in H_{b}\\
(\sigma_{b}^{-1}\alpha_{a}\sigma_{b})(1)=a.
\end{array}\right.\label{eq4}
\end{equation}
By the Lemma \ref{lemma=00003D00002020}, item 2, there exists a
unique permutation $\alpha_{a}$ satisfying (\ref{eq4}). Then the
element $y=\alpha_{a}^{-1}(c)$ is a unique solution of the
equation (\ref{eq=00003D0000203}) in the set $E$, i.e. there exists
a unique solution of the equation (\ref{eq=00003D0000203}) in the
set $E$.

b) Let us have an equation 
\begin{equation}
(x,a,b)=c,\label{eq=00003D0000203-1}
\end{equation}
where $b,c\in E$, $a\in E^{\ast}\setminus\{0,1\}$. Then 
\begin{equation}
\left\{ \begin{array}{c}
c=\alpha_{a}(b)\\
\alpha_{a}\in H_{x}\\
(\sigma_{x}^{-1}\alpha_{a}\sigma_{x})(1)=a.
\end{array}\right.\label{eq4-1}
\end{equation}

Let us assume that there exist $x_{1},x_{2}\in E$, $x_{1}\neq x_{2}$
such that these two elements satisfy (\ref{eq4-1}), i.e. 
\[
\left\{ \begin{array}{c}
\alpha_{a}(b)=c\\
\alpha_{a}\in H_{x_{1}},\quad\alpha_{a}\in H_{x_{2}}\\
(\sigma_{x_{1}}^{-1}\alpha_{a}\sigma_{x_{1}})(1)=a\\
(\sigma_{x_{2}}^{-1}\alpha_{b}\sigma_{x_{2}})(1)=a.
\end{array}\right.
\]
However, 
\[
\alpha_{a}\in H_{x_{1},x_{2}}=St_{x_{1},x_{2}}(G)=<\mathbf{id}>,
\]
because $G$ is a Frobenius group. So $\alpha_{a}=\mathbf{id}$ and we get 
\[
a=(\sigma_{x_{1}}^{-1}\alpha_{a}\sigma_{x_{1}})(1)=(\sigma_{x_{1}}^{-1}\mathbf{id}\sigma_{x_{1}})(1)=\mathbf{id}(1)=1.
\]
The last equality contradicts the condition $a\in E^{\ast}\setminus\{0,1\}$
of this Lemma.

It means that if there exists a solution of the equation (\ref{eq=00003D0000203-1})
for a fixed $a\in E^{\ast}\setminus\{0,1\}$, then it is a unique
solution. Since the set $E$ is finite, then the uniqueness
of the solution of the equation (\ref{eq=00003D0000203-1}) implies
the existence of this solution.
\end{proof}

\section{One-sided \ $S$-system \ of \ operations \ over \ a finite Frobenius group}

According to the results of the previous section, we can construct (over
any finite Frobenius group) a collection of binary operations 
\begin{equation}
A_{a}(x,y)=(x,a,y),\quad\quad a\in E^{\ast}\label{eq=00003D0000207}
\end{equation}
on some finite set $E$, where $E^{\ast}=\left\{ 0,1,...,m\right\} $
is a set of indices of the operations $A_{i}(x,y)$. Moreover, $E^{\ast}\subseteq E$.

\begin{definition}\label{Def2}
A collection of binary operations
$A_{i}(x,y)$ (where $i\in M$, $M$ is a set of indices, $0, \ 1\in M$)
on some set $E$ is called a \textbf{right $S$-system of operations},
if the following conditions hold: 
\begin{enumerate}
\item $A_{0}\left(x,y\right)=x;\quad A_{1}\left(x,y\right)=y;$ 
\item a binary operation $"\circ"$ can be defined on the set $\left\{ A_{i}\right\} _{i\in M}$ of operations $A_{i}$ by the formula 
\[
\left(A_{i}\circ A_{j}\right)\left(x,y\right)=A_{i}\left(i,A_{j}(x,y)\right)=A_{k}\left(x,y\right)\text{ for some }k\in M\setminus\left\{ 0\right\} .
\]
\item the system $\left\langle A_{i}\left(i\neq0\right),\circ,A_{1}\right\rangle $
is a group. 
\end{enumerate}
\end{definition}

\begin{lemma}\label{lemma=00002032}
The collection of binary
operations (\ref{eq=00003D0000207}) is
a right $S$-system of idempotent quasigroups, when $a\neq0,1$.
\end{lemma}

\begin{proof}
Let the conditions of the Lemma hold. Then the item 1 of the Definition \ref{Def2}
follows from (\ref{eq6.8}), and the item 2 of the Definition \ref{Def2} follows
from the item 3 of Lemma \ref{Lemma=00003D00002021}. According to the items
2 and 4 from Lemma \ref{Lemma=00003D00002021} we get that for
$a\neq0,1$ operations $A_{a}(x,y)=(x,a,y)$ are idempotent quasigroups.

According to (\ref{eq=00003D0000203_11}) we have: 
\[
(x,a,(x,b,y))=(x,c,y)
\]
for any $x,y\in E$ and $a,b\in E^{\ast}$, where $c=h_{a}(b),\quad h_{a}=(\sigma_{x}^{-1}\alpha_{a}\sigma_{x})\in H_{0}$.
However $h_{a}h_{b}=h_{c}$ in the subgroup $H_{0}$ of the initial Frobenius
group $G$. Then the system $\left\langle A_{i}\left(i\neq0\right),\circ,A_{1}\right\rangle $
is a group, isomorphic to the group $H_{0}$. It means that the item
3 of the Definition \ref{Def2} holds true.
\end{proof}

\begin{definition}
Two operations $A(x,y)$ and $B(x,y)$on the same
set $E$ are called \textbf{orthogonal}, if, for any $a,b\in E$,
the system 
\[
\left\{ \begin{array}{c}
A(x,y)=a\\
B(x,y)=b
\end{array}\right.
\]
has a unique solution in $E\times E$.
\end{definition}

\begin{lemma}\label{lemma=00002034}(see \cite{Belous65}, p. 104)
For any $a,b\in E^{\ast}$, $a\neq b$, operations $A_{a}(x,y)=(x,a,y)$
and $A_{b}(x,y)=(x,b,y)$ are orthogonal.
\end{lemma}

\begin{proof}
Let $a,b\in E^{\ast}$, $a\neq b$. Let us consider the system 
\begin{equation}
\left\{ \begin{array}{c}
(x,a,y)=c\\
(x,b,y)=d
\end{array}\right.\label{eq=00003D0000206}
\end{equation}
for some arbitrary elements $c,d\in E$. Since the collection of
operations $A_{a}(x,y)=(x,a,y)$ is a right $S$-system, we can represent
the operation $A_{b}(x,y)=(x,b,y)$ in the following way: 
\begin{equation}
(x,b,y)=(x,k,(x,a,y))\label{eq=00003D0000207-1}
\end{equation}
for some parameter $k\in E^{\ast}$.

\emph{Case} $a=0$. This implies $b\neq0$.

If $b=1$ then the system (\ref{eq=00003D0000206}) can be rewritten in
the following form: 
\[
\left\{ \begin{array}{c}
(x,0,y)=c\\
(x,1,y)=d
\end{array}\right.\quad\Leftrightarrow\quad\left\{ \begin{array}{c}
x=c\\
y=d,
\end{array}\right.
\]
i.e. the system (\ref{eq=00003D0000206}) has the unique solution $(c,d)$
in $E\times E$.

If $b\neq0,1$, then the system (\ref{eq=00003D0000206}) can be rewritten
as follows: 
\[
\left\{ \begin{array}{c}
(x,0,y)=c\\
(x,b,y)=d
\end{array}\right.\quad\Leftrightarrow\quad\left\{ \begin{array}{c}
x=c\\
(x,b,y)=d
\end{array}\right.\quad\Leftrightarrow\quad\left\{ \begin{array}{c}
x=c\\
(c,b,y)=d.
\end{array}\right.
\]
Since $b\neq0,1$, the operation $(x,b,y)$ is a quasigroup. Then
the equation $(c,b,y)=d$ has the unique solution $y_{0}\in E$. Then the
system (\ref{eq=00003D0000206}) has the unique solution $(c,y_{0})$
in $E\times E$.

\emph{Case} $a=1$ is proved in a similar manner.

\emph{Case} $a\neq0,1$. If $b=0$, then the system (\ref{eq=00003D0000206})
can be rewritten in the following form: 
\[
\left\{ \begin{array}{c}
(x,a,y)=c\\
(x,0,y)=d
\end{array}\right.\quad\Leftrightarrow\quad\left\{ \begin{array}{c}
(x,a,y)=c\\
x=d,
\end{array}\right.
\]
and the proof is identical to that of the case $a=0,b\neq0,1$.

If $b=1$ then proof is analogical to proof of the previous case.

Let $b\neq0,1$. Then using the formula (\ref{eq=00003D0000207-1}), the
system (\ref{eq=00003D0000206}) can be rewritten as follows: 
\begin{equation}
\left\{ \begin{array}{c}
(x,a,y)=c\\
(x,k,(x,a,y))=d
\end{array}\right.\quad\Leftrightarrow\quad\left\{ \begin{array}{c}
(x,a,y)=c\\
(x,k,c)=d.
\end{array}\right.\label{eq=00003D00002010}
\end{equation}
Since the collection of operations $A_{a}(x,y)=(x,a,y)$ is a right
$S$-system, then the system $\left\langle A_{i}\left(i\neq0\right),\circ,A_{1}\right\rangle $
is a group, isomorphic to the groups $\left\langle E^{\ast}\setminus\left\{ 0\right\} ,\ast,1\right\rangle $
and $H_{0}$. Since $a,b\neq0,1$, $a\neq b$, the formula
(\ref{eq=00003D0000207-1}) implies: 
\[
k\star a=b\quad\Leftrightarrow\quad k=b\star a^{(-1)}\neq0,1,
\]
where $a^{(-1)}$ is an inverse element to the element $a$ in the
group $\left\langle E^{\ast}\setminus\left\{ 0\right\} ,\ast,1\right\rangle $.
Then, by Lemma \ref{Lemma=00003D00002021}, the operation $(x,k,y)$
is a quasigroup. It means that the equation $(x,k,c)=d$ has a unique
solution $x_{0}\in E$. Hence, from (\ref{eq=00003D00002010}), we get
the equation 
\begin{equation}
(x_{0},a,y)=c.\label{eq=00003D0000209}
\end{equation}
As $a\neq0,1$, the operation $(x,a,y)$ is a quasigroup. Since
$x_{0}\in E$ is the determinate fixed element, then the equation (\ref{eq=00003D0000209})
has a unique solution $y_{0}\in E$. So, the pair $(x_{0},y_{0})$ is
a unique solution of the system (\ref{eq=00003D0000206}) in $E\times E$.
\end{proof}

\begin{lemma}Let us define the following two operations $"\star"$
and $"\circ"$ on the set $E^{\ast}$: 
\[
0\star a=a\star0=0,
\]
\[
a,b,c\in E^{\ast}\setminus\left\{ 0\right\} :\quad\quad a\star b=c\quad\overset{def}{\Leftrightarrow}\quad(x,a,(x,b,y))=(x,c,y)
\]
for any $x,y\in E$; and 
\[
a\circ b=c\quad\overset{def}{\Leftrightarrow}\quad(0,a,b)=c.
\]
Then, for any $a,b\in E^{\ast}$ 
\begin{enumerate}
\item $a\star b=a\circ b,$ 
\item the systems $\left\langle E^{\ast}\backslash\left\{ 0\right\} ,\ast,1\right\rangle $
and $\left\langle E^{\ast}\backslash\left\{ 0\right\} ,\circ,1\right\rangle $
are coincided groups. 
\end{enumerate}
\end{lemma}

\begin{proof}
1. We have  
\[
a\star b=c=(0,c,1)=(0,a,(0,b,1))=(0,a,b)=a\circ b
\]
for $x=0,y=1$ and $\forall a,b\in E^{\ast}$ from the definition of operation $"\star"$ and Lemma \ref{Lemma=00003D00002021}.

2. Since collection of operations $A_{a}(x,y)=(x,a,y)$ is a
right $S$-system, then the system $\left\langle A_{i}\left(i\neq0\right),\circ,A_{1}\right\rangle $
is a group, isomorphic to the groups $\left\langle E^{\ast}\backslash\left\{ 0\right\} ,\ast,1\right\rangle $
and $H_{0}$. So the system $\left\langle E^{\ast}\backslash\left\{ 0\right\} ,\ast,1\right\rangle $ coincides with this group.
\end{proof}

Let us remind that $a^{(-1)}$ is an inverse element to the element $a$
in the group $\left\langle E^{\ast}\backslash\left\{ 0\right\} ,\ast,1\right\rangle $.
\begin{lemma}\label{lem34} 
\begin{enumerate}
\item The mappings 
\[
\varphi_{b,a}\left(x\right)=\left(b,a,x\right),\quad b\in E,\quad a\in E^{\ast}\backslash\left\{ 0\right\} ,
\]
\[
\psi_{b,a,d}\left(x\right)=\left(b,a,\left(d,a^{\left(-1\right)},x\right)\right),\quad b,d\in E,\quad a\in E^{\ast}\backslash\left\{ 0\right\} 
\]
are permutations on the set $E.$ 
\item All permutations \ $\varphi_{b,a}$ \ and \ $\psi_{b,a,d}$ \ are permutations from the initial Frobenius group $G$. 
\end{enumerate}
\end{lemma}

\begin{proof}
1. From (\ref{eq6.8},\ref{eq=00003D0000203_11}) we get: 
\[
\varphi_{b,1}\left(x\right)=\left(b,1,x\right)=x,
\]
\[
\psi_{b,a,b}\left(x\right)=\left(b,a,\left(b,a^{\left(-1\right)},x\right)\right)=\left(b,a\ast a^{\left(-1\right)},x\right)=x,
\]
\[
\psi_{b,1,d}\left(x\right)=\left(b,1,\left(d,1,x\right)\right)=x.
\]
If $a\in E^{\ast}\backslash\left\{ 0,1\right\} $, then, for any arbitrary
$b$ and $a$, the mapping $\varphi_{b,a}(x)=L_{b}^{\left(a\right)}\left(x\right)$
is a left translation in the quasigroup $\left(x,a,y\right)$ with
respect to the element $b$, i.e. it is a permutation. Similarly,
for any arbitrary $b,a$ and $d$ we have 
\[
\psi_{b,a,d}\left(x\right)=L_{b}^{\left(a\right)}L_{d}^{\left(a^{(-1)}\right)}\left(x\right),
\]
i.e. the mapping $\psi_{b,a,d}$ is a composition of two translations:
$L_{b}^{\left(a\right)}$ in the quasigroup $\left(x,a,y\right)$
and $L_{d}^{\left(a^{(-1)}\right)}$ in the quasigroup $\left(x,a^{\left(-1\right)},y\right)$. Hence, it is a permutation too.

2. If $a=1$ then $\varphi_{b,1}=\mathbf{id}\in G$. Let $a\in E^{\ast}\backslash\left\{ 0,1\right\} $.
According to the Lemma's conditions and formula (\ref{eq1}) for any
$x\in E$ we have: 
\[
\varphi_{b,a}\left(x\right)=\left(b,a,x\right)=\alpha_{a}(x),
\]
where $\alpha_{a}\in H_{b}\subset G$; i.e. $\varphi_{b,a}=\alpha_{a}\in H_{b}\subset G$.

If $a=1$ or $b=d$, then the item 1 implies: 
\[
\psi_{b,a,b}=\psi_{b,1,d}=\mathbf{id}\in G.
\]
Let $a\in E^{\ast}\backslash\left\{ 0,1\right\} $. Then,for any $b,d\in E,b\neq d$ and for any $x\in E$ we have  
\[
\psi_{b,a,d}\left(x\right)=\left(b,a,\left(d,a^{\left(-1\right)},x\right)\right)=\left(b,a,\varphi_{d,a^{(-1)}}\left(x\right)\right)=\varphi_{b,a}\left(\varphi_{d,a^{(-1)}}\left(x\right)\right),
\]
\[
\psi_{b,a,d}=\varphi_{b,a}\circ\varphi_{d,a^{(-1)}}.
\]
Since $\varphi_{b,a},\varphi_{d,a^{(-1)}}\in G$, then $\psi_{b,a,d}\in G$ too.
\end{proof}

\begin{lemma}\label{lem9.8}
The following statements are true: 
\begin{enumerate}
\item for any $b\in E,\quad a\in E^{\ast}\backslash\left\{ 0,1\right\} $, the
permutation $\varphi_{b,a}$ has one and only one fixed element
$b$, 
\item for any $b,d\in E,\quad a\in E^{\ast}\backslash\left\{ 0,1\right\} $, the
permutation $\psi_{b,a,d}$ is a fixed-point-free permutation,
if $b\neq d$, 
\item for any fixed element $a_{0}\in E^{\ast}\backslash\left\{ 0,1\right\} $
the set of permutations 
\[
T=\left\{ \psi_{b,a_{0},0}|\,b\in E\right\} 
\]
is sharply transitive on the set $E$. 
\end{enumerate}
\end{lemma}

\begin{proof}
1. We have 
\[
\varphi_{b,a}(b)=(b,a,b)=b.
\]
Let us assume that there exists an element $x_{0}\in E$, $x_{0}\neq b$,
such that $\varphi_{b,a}(x_{0})=x_{0}.$Then we obtain 
\[
(b,a,x_{0})=x_{0}.
\]
However, it is clear that 
\[
(x_{0},a,x_{0})=x_{0}.
\]
Since the operation $(x,a,y)$ is a quasigroup and $x_{0}\neq b$,
we have a contradiction between last two equalities. So $\varphi_{b,a}(x)\neq x$
for any $x\in E\backslash\{b\}$.

2. Let us assume that there exists an element $x_{0}\in E$ such
that $b\neq d$ and $\psi_{b,a,d}(x_{0})=x_{0}$, i.e. 
\[
(b,a,(d,a^{(-1)},x_{0})=x_{0}.
\]
Using (\ref{eq=00003D0000203_11}), we get 
\[
\begin{array}{c}
(b,a^{(-1)},x_{0})=(b,a^{(-1)},(b,a,(d,a^{(-1)},x_{0})))=\\
=(b,a^{(-1)}\ast a,(d,a^{(-1)},x_{0}))=(b,1,(d,a^{(-1)},x_{0}))=\\
=(d,a^{(-1)},x_{0}),
\end{array}
\]
i.e. $b=d$, since the operation $(x,a^{(-1)},y)$ is a quasigroup
for any $a\in E^{\ast}\backslash\left\{ 0,1\right\} $. So we have
a contradiction with condition $b\neq d$. It means that if $b\neq d$,
then $\psi_{b,a,d}(x)\neq x$ for any $x\in E$.

3. Let $a_{0}$ be an arbitrary element from $E^{\ast}\backslash\left\{ 0,1\right\} $.
For an arbitrary fixed element $c\in E$, we have: 
\[
\psi_{t,a_{0},0}(c)=(t,a_{0},(0,a_{0}^{(-1)},c))=(t,a_{0},a_{0}^{(-1)}\ast c)=R_{a_{0}^{(-1)}\ast c}^{(a_{0})}(t),
\]
where $R_{b}^{(a)}$ denotes the right translation in the quasigroup
$(x,a,y)$ with respect to the element $b\in E$. As the mapping $R_{a_{0}^{(-1)}\ast c}^{(a_{0})}$
is a permutation on the set $E$, then for any $c,d\in E$ there exists
an unique element $t_{0}\in E$ such that we have 
\[
\psi_{t_{0},a_{0},0}(c)=R_{a_{0}^{(-1)}\ast c}^{(a_{0})}(t_{0})=d,
\]
i.e. the set $T$ is a sharply transitive set of permutations on $E$.
\end{proof}

\begin{lemma}\label{lemma=00002038}
The following statements are true: 
\begin{enumerate}
\item for any $b\in E,\quad a\in E^{\ast}$ 
\[
M_{b}=\left\{ \varphi_{b,a}|\,a\in E^{\ast}\right\} =H_{b},
\] 
\item for any fixed element $a_{0}\in E^{\ast}\backslash\left\{ 0,1\right\} $
the set of permutations 
\[
T_{a_{0}}=\left\{ \psi_{b,a_{0},0}|\,b\in E\backslash\left\{ 0\right\} \right\} 
\]
contains all fixed-point-free permutations from the initial
Frobenius group $G$. Also, $\psi_{0,a_{0},0}=\mathbf{id}.$ 
\end{enumerate}
\end{lemma}

\begin{proof}
1. Using Lemma \ref{lem9.8}, item 1, we get: $\varphi_{b,a}\in H_{b}$
for any fixed $b\in E$. Let us prove that if $a_{1}\neq a_{2}\,\,(a_{1},a_{2}\in E^{\ast})$
then $\varphi_{b,a_{1}}\neq\varphi_{b,a_{2}}$. Let us assume that $a_{1}\neq a_{2}$ but $\varphi_{b,a_{1}}=\varphi_{b,a_{2}}$. It means that 
\[
\alpha_{a_{1}}=\varphi_{b,a_{1}}=\varphi_{b,a_{2}}=\alpha_{a_{2}},
\]
\[
\sigma_{x}^{-1}\alpha_{a_{1}}\sigma_{x}=\sigma_{x}^{-1}\alpha_{a_{2}}\sigma_{x},
\]
\[
h_{a_{1}}=h_{a_{2}},
\]
\[
a_{1}=h_{a_{1}}(1)=h_{a_{2}}(1)=a_{2}.
\]
We get a contradiction. Hence all permutations $\varphi_{b,a_{i}}$
are different. Then 
\[
card\,M_{b}=card\,\left\{ \varphi_{b,a}|\,a\in E^{\ast}\right\} =card\,E^{\ast}=card\,H_{b},
\]
i.e. $M_{b}=H_{b}.$

2. The group $G$ is a finite Frobenius permutation group of degree
$n$, where $n=card\,E$. It is easy to show that the group $G$ contains
exactly $(n-1)$ fixed-point-free permutations (see \cite{Hall62,Wiel64}). It
is true that, for any fixed $a_{0}\in E^{\ast}\backslash\left\{ 0,1\right\} $, the
set $T=\left\{ \psi_{b,a_{0},0}|\,b\in E\right\} $ is transitive
on the set $E$ (Lemma \ref{lem9.8}, item 3), and $\psi_{0,a_{0},0}=\mathbf{id}$
(Lemma \ref{lem34}, item 1). Then the set 
\[
T_{a_{0}}=\left\{ \psi_{b,a_{0},0}|\,b\in E\backslash\left\{ 0\right\} \right\} =T\backslash\left\{ \psi_{0,a_{0},0}\right\} =T\backslash\left\{ \mathbf{id}\right\} 
\]
contains exactly $(n-1)$ fixed-point-free permutations (see Lemma
\ref{lem34}, item 2). So we can see that set $T_{a_{0}}$ contains
all fixed-point-free permutations from the Frobenius group $G$.
\end{proof}

We use a well-known notion of transversal in a group to its subgroup
in next lemma. This notion was described and studied in \cite{Kuzn941}.

\begin{lemma}\label{lemma=00002039}
The following statements are true: 
\begin{enumerate}
\item for any fixed $a_{0}\in E^{\ast}\backslash\left\{ 0,1\right\} $,
the set $T=\left\{ \psi_{b,a_{0},0}|\,b\in E\right\} $, containing
all fixed-point-free permutations from the group $G$ with the identity
element $\mathbf{id}=\psi_{0,a_{0},0}$, is a left transversal in
the group $G$ to subgroup $H_{0}$, 
\item the set $T$ is a loop transversal in $G$ to $H_{0}$, 
\item the set $T$ form a group which is a subgroup of the group $G$, 
\item the subgroup $T$ is a normal subgroup in $G$. 
\end{enumerate}
\end{lemma}

\begin{proof}
1. Let us fix $a_{0}\in E^{\ast}\backslash\left\{ 0,1\right\} $.
Since the set $T=\left\{ \psi_{b,a_{0},0}|\,b\in E\right\} $ is a
transitive set of permutations on the set $E$ (see Lemma \ref{lem9.8},
item 3) then $T$ is a left transversal in $G$ to $H_{0}$.

2. Using Lemma \ref{lem9.8}, item 3, we have that the set $T=\left\{ \psi_{b,a_{0},0}|\,b\in E\right\} $
is a sharply transitive set of permutations on the set $E$. Then $T$
is a left transversal in $G$ to $H_{a}$ for any $a\in E^{\ast}$.
Hence, $T$ is a loop transversal in $G$ to $H_{0}$ (see \cite{Kuzn941},
Lemma 4, item 3 and Lemma 9).

3. From the previous item, we see that the set $T$ is a left transversal
in $G$ to $H_{a}$ for any $a\in E^{\ast}$. Then, for any
$a\in E^{\ast}$ 
\[
t_{i}^{-1}t_{j}\notin H_{a}\quad\forall\;i,j\in E,i\neq j.
\]
It means that for any $i,j\in E,i\neq j$ 
\begin{equation}
(t_{i}^{-1}t_{j})\notin\underset{a\in E^{\ast}}{\bigcup}H_{a}.\label{eq=00003D00002013}
\end{equation}
For any Frobenius group $G$, we have 
\[
G=\left(\underset{a\in E^{\ast}}{\bigcup}H_{a}\right)\bigcup T.
\]
Then (\ref{eq=00003D00002013}) implies that $(t_{i}^{-1}t_{j})$
is a fixed-point-free permutation, i.e. 
\[
t_{i}^{-1}t_{j}=t_{k}
\]
for some element $t_{k}\in T$, since the set $T$ contains all the
fixed-point-free permutations of the group $G$. This means that the set
$T$ is closed under group multiplication (and inverse transformation) in the group $G$. Hence, the set $T$ form a group which is a subgroup of $G$.

4. This statement is obvious, because any permutation conjugated to a fixed-point-free permutation is also a fixed-point-free permutation.
\end{proof}

This completes the character-free proof of the Frobenius theorem.

\end{document}